\documentstyle[amscd,amssymb,verbatim,11pt]{amsart}

\theoremstyle{plain}
\newtheorem{Prop}{Proposition}[section]
\newtheorem{Thm}[Prop]{Theorem}
\newtheorem{Cor}[Prop]{Corollary}

\newtheorem{Lem}[Prop]{Lemma}

\newtheorem{Qs}[Prop]{Question}

\theoremstyle{definition}
\newtheorem{Def}[Prop]{Definition}

\theoremstyle{remark}
\newtheorem{Rem}[Prop]{Remark}

\newtheorem{Exam}[Prop]{Example}

\def\int{\mathop{\roman{int}}}

\def\1{^{-1}}

\def\int{\text{Int}}
\errorcontextlines=0 \numberwithin{equation}{section}
\begin{document}
\title[
Cohomological dimension with respect to nonabelian groups
]%
   {Cohomological dimension with respect to nonabelian groups
}

\author{Atish Mitra}
\address{University of Tennessee, Knoxville, TN 37996, USA}
\email{ajmitra@@math.utk.edu}

\date{ June 28, 2005
} 


\thanks{
}

\begin{abstract}
Cencelj and Dranishnikov showed \cite{CD1,CD2,CD3}that for certain
nilpotent groups $G$,  $K(G_{ab},1) \in \text{AE}(X)$ is equivalent
to $K(G,1) \in \text{AE}(X)$ for any compacta $X$ (here $G_{ab}$ is
the abelianization of $G$). We examine the same problem  for
solvable groups. We also give an elementary proof of this fact for
any nilpotent group and any 2-dimensional metric space.
\end{abstract}

\maketitle

\medskip
\medskip
\tableofcontents

\section{Introduction}

We define cohomological dimension $\text{dim}_G X$ of a topological
space $X$ with respect to an abelian group $G$ as the largest number
$n$ such that there exists a closed subset $A \subset X$ with
$\check{H}^n(X,A;G) \ne 0$. Several authors have worked on
developing this theory for various classes of spaces like compacta
and metrizable spaces. For surveys of the theory the reader is
referred to \cite{D1}, \cite{Dy1}.

We use Kuratowski's  notation $X \tau M$ if every map from a closed
subset of $X$ to $M$ extends over $X$. The following theorem allows
us to formulate an equivalent definition of cohomological dimension
of compacta with respect to abelian groups in terms of extension of
maps to Eilenberg-Maclane  spaces.

\begin{Thm}\cite{D2}
For any compactum $X$ and abelian group $G$ the following conditions
are equivalent:
\begin{enumerate}
\item $\text{dim}_G X \le n$
\item $\check{H}^{n+1}(X,A;G) = 0$ for all closed $A \subset X$
\item $X \tau \text{K}(G,n)$
\end{enumerate}

\end{Thm}

Trying to extend the notion of cohomological dimension to
non-abelian groups $G$, one notices that the corresponding
Eilenberg-Maclane spaces are defined only for $n=1$.Dranishnikov and
Repovs considered the case of perfect groups in \cite{DR}. Later
Cencelj and Dranishnikov studied this notion in a series of three
papers \cite {CD1,CD2,CD3}.

Generalizing Dranishnikov's theorem about extension of maps to
simply-connected complexes \cite{D1}, they obtained the following
result.

\begin{Thm}
For any nilpotent CW-complex $M$ and finite-dimensional metric
compactum $X$, the following are equivalent:
\begin{enumerate}
\item $X \tau M$
\item $X \tau {\text{SP}}^{\infty}M$
\item ${\text{dim}}_{H_i(M)} X \le i$ for every $i > 0$
\item ${\text{dim}}_{{\pi}_i(M)} X \le i$ for every $i > 0$
\end{enumerate}
\end{Thm}

Recall that ${\text{SP}}^{i} X = {X^i}/({\Sigma}_i)$, where
${\Sigma}_i$ is the symmetric group on $i$ letters. If $X$ is
pointed, there is a natural embedding ${\text{SP}}^{i} X
 \to {\text{SP}}^{i+1} X$, and ${\text{SP}}^{\infty} X $ is defined
 as direct limit of the ${\text{SP}}^{i} X$.

In course of the proof of the above theorem, they obtained the
following theorem.

\begin{Thm}
For a nilpotent group $N$ and every metric compactum $X$ the
following equivalence holds : $\text{dim}_N X \le 1 \Leftrightarrow
{\text{dim}_{N_{ab}}} X \le 1$ provided $N $ has one of the
following properties:

\begin{enumerate}

\item $N$ is a torsion group.
\item for every prime $p$ s.t. $\text{Tor}_p N \neq 1$
\begin{enumerate}
\item $N$ is not $p$-divisible, or
\item $\text{Tor}_p {N_{ab}} \neq 0$.
\end{enumerate}
\end{enumerate}

\end{Thm}

In view of Theorem 1.3 it is  worthwhile studying the following
problem.

\begin{Qs}Let $\mathcal C$ be a class of spaces. Describe groups for which the
following are equivalent for any $X \in \mathcal C$?
\begin{enumerate}
\item $X \tau K(G,1)$
\item $X \tau K(G_{ab},1)$
\end{enumerate}
\end{Qs}

In section 2 we list well known results which we use in the sequel.
In section 3 we show that if $\mathcal C$ is the class of all
2-dimensional metrizable spaces, then nilpotent groups satisfy
question 1. In section 3 we show that if $\mathcal C$ is the class
of all 2-dimensional compacta, then finite solvable groups satisfy
question 1.4. Finally we  give an example showing that infinite
solvable groups do not satisfy question 1.4 when if $\mathcal C$ is
the class of all 2-dimensional compacta.
\section{Preliminaries}

\vspace{.2 in} We collect below some classical results which we will
use in the sequel.

\begin{Prop}
Let $A$ be closed in $X$. A map $f:A \to Y$ has an extension to $X$
iff some homotopic map $g:A \to Y$ has an extension.
\end{Prop}
\vspace{.1 in}

Recall that a map $f: X \to Y$ between CW complexes is called
cellular if $f(X^{(n)}) \subset Y^{(n)}$ for all $n$.
\begin{Thm}[\textbf{Cellular Approximation Theorem}]\cite{LW}
Every map $f: X \to Y$ between CW complexes is homotopic to a
cellular map.
\end{Thm}


\begin{Thm}[\textbf{Seifert - Van Kampen }]
Let $X_1$, $X_2$ be path connected open subsets of $X$,  and the
intersection $X_0=X_1 \cap X_2$ be path connected. Then the
following commutative diagram of fundamental groups (based at some
$x_0 \in X_0$) where all maps are induced by inclusions is a pushout
diagram:

\[
\begin{CD}
     {\pi}_1 (X_0)   @>{i_1}>> {\pi}_1 (X_1) \\
     @V{i_2}VV                  @VV{j_1}V \\
     {\pi}_1 (X_2)   @>{j_2}>> {\pi}_1 (X)
\end{CD}
\]

\end{Thm}


\begin{Thm}[\textbf{Bockstein exact sequence}]\cite{Dy1}
If $1 \to G \to E \to \Pi \to 1$ is a short exact sequence of
abelian groups, then there is a natural exact sequence

\begin{equation}
   \cdots \to  \check{H}^n(X,A;G) \to \check{H}^n(X,A;E) \to
   \check{H}^n(X,A;\Pi) \to
   \check{H}^{n+1}(X,A;G) \to \cdots
\end{equation}

for any paracompact space $X$ and its closed subspace A. Here $
\check{H}^n(X,A;G)$ is the $n$-th Cech cohomology group of the pair
$(X,A)$ with coefficients in $G$.

\end{Thm}


The following theorem due to Cenceli and Dranishnikov  gives one of
the implications of question 1 for any group $G$

\begin{Thm}\cite{CD1}
Let $X$ be a compactum and $G$ a group. Then $X \tau K(G,1)
\Rightarrow X \tau K(G_{ab},1)$.
\end{Thm}

\emph{Proof:} From \cite{D1}, theorem 6 it follows that $X \tau
{\text{SP}}^{\infty}K(G,1)$. We use esults from \cite{DT}, \cite{H}
to see that $X \tau K(H_1(K(G,1)),1)$, which is equivalent to $X
\tau K(G_{ab},1)$.


\section{Nilpotent Groups}

We present a generalization of theorem 1.3 for 2-dimensional
metrizable spaces.

\begin{Thm}
Let $G$ be a nilpotent group and $X$ be a 2-dimensional metric
space.
\\
Then $X \tau K(G,1) \Leftrightarrow X \tau K(G_{ab},1)$.
\end{Thm}

Recall that the lower central series of a group $G$ is $G={\gamma}_1
G \geqq {\gamma}_2 G \geqq \cdots $, where ${\gamma}_{i+1} G =
[{\gamma}_{i} G, g]$. Notice that ${\gamma}_i G / {\gamma}_{i+1} G$
lies in the center of $G / {\gamma}_{i+1} G$ and that the nilpotent
class of $G$ equals the length of the lower central series of $G$.
The following lemma shows that the first lower central factor
$G_{ab}$ exerts a strong influence on the subsequent  lower central
factors.

\begin{Lem}[Robinson]

Let G be a group and let $F_i={\gamma}_i G / {\gamma}_{i+1} G$. Then
the mapping $a({\gamma}_{i+1} G) \otimes g[G,G] \mapsto
[a,g]({\gamma}_{i+2} G)$ is a well defined module epimorphism from
$F_i \otimes _{{\mathbb {Z}}} G_{ab}$ to $F_{i+1}$.

\end{Lem}

Proof: Let $g \in G$ and $a \in G_i$, and consider the function
$(a({\gamma}_{i+1} G) , g[G,G]) \mapsto [a,g]({\gamma}_{i+2} G)$.
Check that the mapping is well-defined and bilinear. By the
universal property of tensor products, we have an induced
homomorphism. It is easily checked that it is an epimorphism.

\vspace{.2 in}

The above lemma gives the following useful corollary.

\begin{Cor}

The ${\cal{P}}$ be a group theoretical property which is inherited
by images of tensor products (of abelian groups) and by extensions.
If $G$ is a nilpotent group such that $G_{ab}$ has ${\cal{P}}$, then
$G$ has ${\cal{P}}$.

\end{Cor}

Proof: Let $F_i={\gamma}_i G / {\gamma}_{i+1} G$, and note that
$F_1=G_{ab}$. Suppose $F_i$ has ${\cal{P}}$ : then the previous
lemma implies that $F_{i+1}$ has ${\cal{P}}$. As $G$ is nilpotent,
some ${\gamma}_{c+1}=1$. Since ${\cal{P}}$ is closed under
extensions, $G$ has ${\cal{P}}$.

\vspace{.2 in}

\emph{Proof of theorem 3.1:} Let ${\cal{P}}$ be the  group
theoretical property of "belonging to ${\cal{G}}$". We know that if
$1 \to G \to E \to \Pi \to 1$ is a short exact sequence and $X \tau
K(G,1)$ and $X \tau K(\Pi,1)$, then $X \tau K(E,1) $\cite{Dy1}. We
also know that if $X$ is metrizable, then $X \tau K(G,1) \Rightarrow
X \tau K(G \otimes H,1)$ for any abelian groups $G,H$ \cite{Dy1}.
From the Bockstein exact sequence we see that if $X$ is
2-dimensional, $G$ abelian and $X \tau K(G,1)$, then $X \tau
\phi(K(G,1))$ for any homomorphism $\phi$. Apply the previous
corollary.

\begin{Rem}
The theorem can be improved to the case where $X$  has cohomological
dimension 2.

\end{Rem}

\section{Finite Solvable Groups}

Our effort to answer question 1 for solvable groups gives the
following theorem.

\begin{Thm}
Let $G$  be a finite solvable group and $X$ be a 2-dimensional
metrizable space. Then $X \tau K(G,1) \Leftrightarrow X \tau
K(G_{ab},1)$.
\end{Thm}

For proving the above theorem, we need some lemmas.
\begin{Prop}
Let $X$ be a metrizable space and $K_1$, $K_2$ be subcomplexes of
the CW complex $K_1 \cup K_2$. Then:
\begin{enumerate}
\item If $K_1$, $K_2$ and $K_0=K_1 \cap K_2$ are absolute extensors
of $X$, then so is $K_1 \cup K_2$.
\item If $K_1 \cup K_2$ and $K_0=K_1 \cap K_2$ are absolute extensors
of $X$, then so are $K_1$, $K_2$.
\end{enumerate}
\end{Prop}

\emph{Proof:} (1) Let $f: A \to K_1 \cup K_2$ be a map from a closed
subset $A$ of $X$. Define $C_i=f^{-1}(K_i)$ and note that $A=C_1
\cup C_2$ and $f(C_1 \cap C_2) \subset {K_1 \cap K_2}$. As the
closure of any one of the sets $C_1-C_2$ and $C_2-C_1$ misses the
other (and as $X$ is completely normal), we can find  $U$ open in
$X$ such that $C_1-C_2 \subset U \subset \overline{U} \subset X -
(C_2-C_1)$. Define $D_1=\overline{U} \cup (C_1 \cap C_2)$ and
$D_2=(X-U) \cup (C_1 \cap C_2)$. Then we have $D_1 \cap A = C_1$,
$D_2 \cap A = C_2$ and $D_1 \cup D_2= X$. As $C_1 \cap C_2$ is
closed in $D_1 \cap D_2$ (which is closed in $X$) and as $K_1 \cap
K_2$ is absolute extensor of $X$, we can extend $f \mid C_1 \cap
C_2$ to $\overline{f}:D_1 \cap D_2 \to K_1 \cap K_2$. Now we have a
map from $C_1 \cup (D_1 \cap D_2)$ to $K_1$ and a map from $C_2 \cup
(D_1 \cap D_2)$ to $K_2$, both of which we can extend and then paste
the resulting maps to get the desired extension.

(2) Start with $f: A \to K_1$ and extend to $\overline{f}: X \to K_1
\cup K_2$. Calling $C_i=(\overline{f})^{-1}(K_i)$, we have $X=C_1
\cup C_2$, $A \subset C_1$ and $\overline{f} (C_1 \cap C_2) \subset
K_1 \cap K_2$. Extend $\overline{f} \mid C_1 \cap C_2$ to
$\widetilde{f}: C_2 \to K_1 \cap K_2$. Paste $\overline{f} \mid C_1$
and $\widetilde{f}\mid C_2$ to get the desired extension.

\begin{Def}
Let $G_1$, $G_2$ and $A$ be groups, and $f_i:A \to G_i$ be
homomorphisms. Then the amalgamated product $G_1 *_{A}G_2$ is
defined as $\displaystyle{\frac{G_1 * G_2}{  \langle f_1(a)=f_2(a),
 a \in A \rangle}}$(Here $\langle \rangle$ denotes the normal
closure).  Notice that the amalgamated product is a pushout of the
following diagram:

\[
\begin{CD}
     A   @>{f_1}>> G_1 \\
     @V{f_2}VV                   \\
     G_2
\end{CD}
\]

\end{Def}
\vspace{.2 in}

\begin{Prop}
Let $X$ be a 2-dimensional metrizable space and define ${\cal{H}}_X$
$= \{G:X\tau K(G,1)\}$.
\begin{enumerate}
\item If $G_1$,$G_2$ and $A$ are in ${\cal{H}}_X$ , then so is $G_1 *_{A}G_2$.
\item If $G_1 *_{A}G_2$ and $A$ are in ${\cal{H}}_X$, then so are $G_1$ and $G_2$.
\end{enumerate}
\end{Prop}

\emph{Proof:} Note that a homomorphism $f:A \to G$ induces a map
$f_{*}:K(A,1)^{(2)} \to K(G,1)^{(2)}$: we start by sending the
generators to the corresponding loops in  $K(G,1)^{(1)}$, and then
use the relators of $A$ to extend the map to the 2-skeleton . Let
$M_{f_{*}}$ be the mapping cylinder of $f_{*}$: then we have
$K(A,1)^{(2)} \hookrightarrow M_{f_{*}}$. Apply the previous theorem
with $K_1=M_{f_{*}1}$, $K_2=M_{f_{*}2}$ and $K_0=K(A,1)^{(2)}$.

\vspace{.5 in }
 We can now prove the following lemma, which we use
later to see that any finite group with solvability index 2 is in
$\cal{G}$.

\begin{Lem}
Let $X$ be a 2-dimensional metrizable space with $\text{dim
}_{{\mathbb Z }/ p} X =1$ and let $K$ be a CW complex. Let $\alpha
\in {\pi}_1(K)$ with ${\alpha}^{p}=1$, for some prime $p$. Then $X
\tau K \Leftrightarrow X \tau K \cup_{\alpha} D^2$, where $D^2$ is a
2-cell.
\end{Lem}

\emph{Proof:} Note that $\alpha:S^1 \to K$ induces a homomorphism
${\alpha}_{\ast}: {\mathbb{Z}} \to {\pi}_1 (K)$ with
${{\alpha}_{\ast}(1)}^p=1$, i.e. ${\alpha}_{\ast}(p)=1$. Then
${\pi}_1(K \cup_{\alpha} D^2)$ is the pushout of

\[
\begin{CD}
     {\mathbb{Z}}   @>{\alpha_1}>> {\pi}_1(K) \\
     @VVV                   \\
     \{1\}
\end{CD}
\]

But in this case it is also the pushout of

\[
\begin{CD}
     {\mathbb{Z}}/p   @>{\alpha_1 \mid}>> {\pi}_1(K) \\
     @VVV                   \\
     \{1\}
\end{CD}
\]

Then we can use the previous theorem.

\vspace{.5 in}

The previous lemma has the following generalization:

\begin{Lem}
Let $X$ be a 2-dimensional metrizable space , ${\cal{P}} = \{p :
\text{dim }_{{\mathbb Z }/ p} X =1\}$ and let $K$ be a CW complex.
Let ${\alpha}_i \in {\pi}_1(K)$ with ${{\alpha}_i}^{{p_i}^{m_i}}=1$,
for some prime $p_i \in {\cal{P}}$ and positive integer $m_i$. Then
$X \tau K \Leftrightarrow X \tau K \cup_{{\alpha}_i} {D_i}^2$, where
${D_i}^2$ is a 2-cell.
\end{Lem}

\emph{Proof of theorem 4.1:}

$A=G_{ab}$ is  a finite abelian group and hence of the form
${\mathbb Z}/{{q_1}^{n_1}} \oplus \cdots \oplus {\mathbb
Z}/{{q_l}^{n_l}}$. Using the previous lemma, we can "kill" all
elements of order products of $q_i$: we create the quotient of $G$
by the normal subgroups generated by those elements. As $G$ is
finite , we can repeat this process a finitely many times to get a
solvable group $\tilde{G}$ such that $X \tau K(G,1) \Leftrightarrow
X \tau K(\tilde{G},1)$ for 2-dimensional compacta $X$.

We claim that $\tilde{G}$ has  trivial abelianization: in the exact
sequence $1 \to [\tilde{G},\tilde{G}] \to \tilde{G} \to
{\tilde{G}}_{ab} \to 1$ let $y$ be any generator of
${\tilde{G}}_{ab}$ and consider it's preimage $y^{'}$. As $y^q=1$
for some $q \in \{q^1, \cdots, q^l\}$, ${y^{'}}^q$ is mapped to 1 in
${\tilde{G}}_{ab}$, and hence  ${y^{'}}^q=k$ for some element $k \in
[\tilde{G},\tilde{G}]$. Every element in $[\tilde{G},\tilde{G}]$ has
order relatively prime to the $q_i$(s), so suppose $k$ has order
$p$, relatively prime to the $q_i$(s). Thus we have ${y^{'}}^{pq}=1$
which implies ${y^{'}}^p=1$, which in turn implies $y^p=1$. Thus
$y=1$, and our claim is proved.

\begin{Rem}
The author thinks that the above theorem is true for any solvable
group with all elements of finite order. There are a few details
which are being worked on.  This paper is a preliminary version
before the "22nd annual workshop in geometric topology, Colorado
College" and the Bedlewo conference.

\end{Rem}

\section{Extension Properties of the Pontryagin Disk}

We will use for the sequel a 2-dimensional compactum called the
mod-2  Pontryagin disk ${{\mathbb{P}}D}^2$, which is a variant of a
classical construction of L.S. Pontryagin \cite{P}.

The space ${{\mathbb{P}}D}^2$  is constructed as the inverse limit
of an inverse sequence of spaces. As the first stage of the
construction, we start with the usual 2-cell embedded in ${\mathbb
R}^3$ and take some triangulation of it with mesh less than 1. To
get the $k+1$th stage from the  $k$th stage we take each 2-simplex,
remove it's interior and glue a disk along the degree 2 map, then
give a triangulation of mesh less than $1/(k+1)$ to the resulting
space. For the bonding map from  $k+1$th stage to the  $k$th stage,
we send the interior 2-simplices of the attached disks to the
barycenter of the removed 2-simplices.

The following theorem is very useful in constructing compacta with
special extension properties. For explanation of terms and proof see
\cite{DR}. Note that $\text{mesh}\{{\lambda}_{i}\} \to 0$ means that
for every $k$, \\ ${\lim}_{i \to \infty}
\{\text{mesh}(q^{k+i}_k({\lambda}_{k+i})\}=0$, which is not the same
as ${\lim}_{i \to \infty} \text{mesh} {\lambda}_i=0$.

\begin{Thm}
Suppose that $K$ is a  countable CW complex and that $X$ is a
compactum such that $X=\lim_{\leftarrow} (X_i,f_{i}^{i+1})$, where
$(L_i,f_{i}^{i+1})$ is a $K$-resolvable inverse system of compact
polyhedra $L_i$ with triangulations ${\lambda}_i$ such that
$\text{mesh}\{{\lambda}_i\} \to 0$. Then $X \tau K$.
\end{Thm}




\begin{Thm}
Let  ${{\mathbb{P}}D}^2$ be the  Pontryagin disk. Then
${{\mathbb{P}}D}^2 \tau K({\mathbb{Z}}/2,1)$.

\end{Thm}

Proof: Let $A \subset {{\mathbb{P}}D}^2$ be a closed and consider
the map $f: A \to K({\mathbb{Z}}/2,1)$. Represent
$({{\mathbb{P}}D}^2,A) = \lim_{\leftarrow}
\{(({{\mathbb{P}}D}^2)_i,A_i),q_{i}^{i+1}\}$  in the usual way as
the inverse limit of polyhedral pairs (all lying in $I^{\infty}$ ).
Like in the preceding theorem, find $i_0$ such that $f$ can be
extended over $A_{i_0}$ to give a map $\overline{f_{i_0}}:A_{i_0}
\to K({\mathbb{Z}}/2,1)$ with $\overline{f_{i_0}} \circ
(q_{i_0}^{\infty}) \simeq f$. Extend $\overline{f_{i_0}} $ over the
1-skeleton of $(({{\mathbb{P}}D}^2)_i$ and get an induced extension
upto homotopy  $\overline{G}: A \cup  ({{\mathbb{P}}D}^2)^{(1)} \to
K({\mathbb{Z}}/2,1)$. Pick any 2-simplex ${\sigma}_{j} \in
(({{\mathbb{P}}D}^2)_{i_0})^{(1)} - (A_{i_0})^{(1)}$ and consider
${\gamma}_{{\sigma}_j} =g_{i_0}(\partial {\sigma}_j) \in
{\mathbb{Z}}/2$. Consider $(q_{i_0}^{\infty})^{-1} ({\sigma}_j)
\subset {{\mathbb{P}}D}^2$ and note that we have an extension (of
$f$) over $A \cup (q_{i_0}^{\infty})^{-1} (\partial {\sigma}_j)$. We
want to extend over $A \cup (q_{i_0}^{\infty})^{-1} ( {\sigma}_j)$:
if ${\gamma}_{{\sigma}_j}$ is trivial there is no obstruction. If
not, find the smallest index $i_1 > i_0$ such that
$(p_{i_0}^{i_1})^{-1}$ has an $({\mathbb{R}}P)^2$ (needs explanation
/ rephrasing ). Clearly the extension over $A \cup
(q_{i_0}^{\infty})^{-1} $ is possible. Do the same for all such
2-simplices ${\sigma}_j)$.

\section{Infinite Solvable Groups}

\begin{Thm}

Let  ${{\mathbb{P}}D}^2$ be the  Pontryagin disk and $\Gamma$ be a
torsion free  group with abelianization being a non-trivial 2-group.
Then $ K({\Gamma}_{ab},1)$ is an absolute extensor of
${{\mathbb{P}}D}^2$, but $ K(\Gamma,1)$ is not an absolute extensor
of ${{\mathbb{P}}D}^2$.
\end{Thm}

\emph{Proof:} Let ${\sigma}$ be any 2-simplex in the first stage of
construction of ${{\mathbb{P}}D}^2$ and identity $\partial \sigma $
and its image under in $X$ under $(f^{\infty}_0)^{-1}\partial \sigma
$. Consider the map $f:\partial {\sigma} \to K(\Gamma,1)$ that sends
$\partial {\sigma}$ to the generator $x$  along the identity map
$S^1 \to S^1$.

We claim that this map does not extend: suppose there is an
extension $g:{{\mathbb{P}}D}^2 \to K(\Gamma,1)$. Then there is some
integer $i_0$ such that there is an extension up to homotopy
$g_{i_0}:{{{\mathbb{P}}D}^2}_{i_0} \to K(\Gamma,1)$. In that case,
$x$ could be expressed as products of elements of $\Gamma$, each of
order 2. But $\Gamma$ is torsion-free, and we have a contradiction.


We give 2 examples of groups satisfying theorem 6.1, with
${\Gamma}_{ab}$ being a 2-group.

\begin{Exam}
Consider the following group, which was first studied by Hirsch:\\

${\Gamma}_1 = \langle x,y,z | x^z=x^{-1},y^z=y^{-1}, [x,y]=z^4
\rangle.$\\

This example of Hirsch is a torsion free polycyclic group which is
not poly infinite cyclic. The  folowing lemma shows that $G$ is not
nilpotent.
\end{Exam}

\begin{Lem}[Robinson, 5.2.20]
A finitely generated torsion-free nilpotent group has a central
series with infinite cyclic factors.
\end{Lem}

The author thanks Prof. Thomas Farrell for pointing out the next
example.

\begin{Exam} It is known \cite{W}{theorem 3.5.5} that there are just 6 affine diffeomorphism classes of compact
connected orientable flat 3-dimensional Riemannian manifolds,
represented by manifolds ${\mathbb R}^3/ G$, where $G$ is one of 6
given groups. One of them is our next example\\

${\Gamma}_2 = \langle \alpha, \beta, \gamma x,y,z | \gamma \beta
\alpha  =xy,  {\alpha}^2=x, \alpha y {\alpha}^{-1}=y^{-1}, \alpha z
{\alpha}^{-1}=z^{-1},\\
 \beta x {\beta}^{-1}=x^{-1},  {\beta}^2=y, \beta z
{\beta}^{-1}=z^{-1},\gamma x {\gamma}^{-1}=x^{-1},\\
\gamma y {\gamma}^{-1}=y^{-1},{\gamma}^2=z \rangle.$\\

It's abelianization is ${\mathbb Z}/4 \oplus {\mathbb Z}/4$. It is
torsion free and  solvable but not nilpotent.

\end{Exam}

\begin{Rem}Even solvable groups with some elements of infinite order may not be
enough to serve as an example of groups in theorem 6.1.\\ As an
example, consider $D_{\infty}$, which is defined to be the
semi-direct product ${\mathbb{Z}}\rtimes_{\theta} {\mathbb{Z}}/2$,
where $\theta:{\mathbb{Z}}/2 \to \text{Aut}(Z)$ is the "$t^{-1}$"
map. A presentation of $D_{\infty}$ is $\langle x,y \mid y^2=1,
yxy^{-1}=x^{-1}\rangle$. Note that the relations imply $(yx)^2=1$,
so $D_{\infty} ={\mathbb{Z}}/2 \ast {\mathbb{Z}}/2$ and  therefore,
the abelianization is ${\mathbb{Z}}/2 \oplus {\mathbb{Z}}/2$. The
commutator is $\langle x \rangle = {\mathbb{Z}}$.

As $D_{\infty}={\mathbb{Z}}/2 \ast {\mathbb{Z}}/2$,
${{\mathbb{P}}D}^2 \tau K(D_{\infty},1)$.

\end{Rem}

\enddocument